\documentclass[journal,12pt,onecolumn,draftclsnofoot]{IEEEtran}
\usepackage{amssymb}
\usepackage{tikz,url}

\newtheorem{thm}{Theorem}

\newtheorem{cor}{Corollary}
\newtheorem{defn}{Definition}

\newcommand{\A}{\mathcal{A}}

\ifCLASSINFOpdf
\else
\fi
\usepackage{amsmath}
\usepackage{multirow}

\usepackage{fixltx2e}

\usepackage{stfloats}
\begin{document}
%
\title{On the Weights of General MDS Codes}
%
%
%

\author{Tim~Alderson~\IEEEmembership{}
\thanks{T. Alderson is with the Department
of Mathematics and Statistics, University of New Brunswick Saint John, Saint John,
NB, E2L 4L5 Canada e-mail: Tim@unb.ca.}
\thanks{Copyright (c) 2017 IEEE. Personal use of this material is permitted.  However, permission to use this material for any other purposes must be obtained from the IEEE by sending a request to pubs-permissions@ieee.org.}
}

%
%

\markboth{Manuscript 2019}%
{ : Bare Demo of IEEEtran.cls for IEEE Journals}
%



\maketitle

\begin{abstract}
The weight spectra of MDS codes of length $ n $ and dimension
$ k $ over the arbitrary alphabets are studied.  For all $ q $-ary MDS codes of dimension $ k $ containing the zero codeword, it is shown that all $ k $ weights from $ n $ to $ n-k+1 $ are realized. The remaining case $ n=q+k-1 $ is also determined.  Additionally, we prove that all binary MDS codes are equivalent to linear MDS codes.  The proofs are combinatorial, and self contained.
\end{abstract}


%
\IEEEpeerreviewmaketitle

\section{Introduction}
%
%
%
%
\IEEEPARstart{T}{he}   weight spectrum of an error correcting code provides the non-zero weight parameters among all codewords, whereas the weight distribution counts, for every given weight parameter, the number of codewords of that weight.  Discussions on weight spectra of codes  can be traced to \cite{MacWilliams1963}, where the author provides partial answers to  the following question. Given a set of positive integers, $ S $,  is it possible to construct a code whose set of non-zero weights is $ S $? 
More recent works, such as \cite{Shi2018},\cite{Shi2018a}, \cite{AldFWS},  and \cite{Alderson2019} investigate upper bounds on the size of the weight spectra of linear codes, and the existence of  maximum weight spectrum (MWS) codes, and full weight spectrum (FWS) codes.   Shi et. al. \cite{Shi2018} skirmish with the analogous problem for non-linear codes, managing to provide an asymptotic bound. Determining general weight spectra for nonlinear codes seems less tractable.

The weight spectrum/distribution of a code are key characteristics  which are of both theoretical and practical interest. The error correcting capability of a code is most typically  discussed in relation to it's minimum distance relative to it's length, as this directly relates to the number of independent errors that can be corrected. However, the weight distribution  and/or  weight spectrum of a code can also provide insight into the error correcting properties of a code. For example, many algorithms for list-decoding such as those in \cite{Goldreich1989,Goldreich2000,Sudan1999}  rely extensively in in determining and bounding the non-empty list sizes  for different distance parameters. The weight spectrum of a code can also provide insights regarding connections with designs \cite{Assmus1969}, and with regard to the extendability of the code \cite{Hill1999,Hill1995}. \\
 
 In \cite{Ezerman2011}, Ezerman, Grasl, and Sol\'{e}  determine the weight spectra of certain linear MDS codes, namely those that satisfy the MDS Conjecture. A key motivation for their investigation is a construction of R\"{o}tteler et. al. \cite{Rotteler2004} \cite{GRASSL_2004} providing the existence of quantum MDS codes for all lengths $ n \le q+1 $.  The construction hinges on the work of Rains \cite{Rains1999}, which in turn relies on the following assertion: 
\begin{itemize}
	\item [($ \star $)]\textit{An $ [n,k,d]_q $-MDS code  has precisely $ k $ distinct nonzero weights, $ n,n-1,\ldots,n-k+1$}. 
\end{itemize}

%

This assertion seems reasonable enough, and may be found for example in [\cite{MR0465509}, p. 320] without an explicit proof. Ezerman et. al \cite{Ezerman2011}   determined that among linear MDS codes  there were in fact some exceptions to ($ \star $), these being the dual of the binary repetition code of length $ n > 2 $ which contains only words of even weight; the $ q $-ary simplex code with parameters $ [q + 1, 2, q]_q $ which contains only words of weight zero or $ q $, and the $ [q+2,3,q]_q $ codes where $ q $ is necessarily even, having two nonzero weights, $ q $, and $ q-2 $. This served to validate the QECC constructions in \cite{Rotteler2004} \cite{GRASSL_2004}, as they did not employ these exceptional codes.\\
In this communication we determine the weight spectrum of all MDS codes, linear or not. In particular, our results subsume those of \cite{Ezerman2011},  and determine all parameters for which ($ \star $)  holds.   


\section{Preliminaries}

Let $ \A $ be an alphabet of size $q$, where $ q $ need not be a prime power. The \textit{Hamming distance}, $ d(a,b) $, between two elements $a,b \in \A^n$ is the number of coordinates in which $ a $ and $ b $ differ. The Hamming distance induces a metric on $\A^n$.  A subset $ C\subseteq \A^n $ with $ |C|=q^k $ is an $(n,k)_q$ \emph{code} when equipped with the Hamming distance, the elements of $ C $ are called \textit{codewords}.  A code is \emph{linear} if its codewords form a $k$-dimensional subspace of the vector space $\mathbb{F}_q^n$ where $\mathbb{F}_q$ denotes the field of order $q$. Unless otherwise stated, the codes discussed in the sequel are not assumed to be linear.\\
The \emph{minimum distance} $d$ of $C$ is the minimum over all distances between distinct pairs of codewords. We denote such a code by an $ (n,k,d)_q $-code.  From the Singleton bound, we obtain
$ d\le n-k+1$. If $C$ meets the Singleton bound, it is said to be \emph{Maximum Distance Separable} (MDS), denoted $(n, k)_q$-MDS (or $[n, k]_q$-MDS in the linear case). The main conjecture on MDS codes states that if $ k \le q $ then a linear MDS code has length
$ n \le q + 1 $ unless $ q $ is even  and $ k = 3 $ or $ k = q - 1 $, in which case $ n \le q + 2 $. The conjecture can be stated without reference to MDS codes, and was first considered by Beniamino Segre in the 1950s \cite{MR0075608}. We remark that there are no known MDS codes (linear or not) which violate the MDS conjecture. In fact, the author would go so far as to propose that \textit{the MDS conjecture holds for non-linear MDS codes as well}.

In what follows we shall require the notion of equivalence among codes. Two codes are said to be \emph{equivalent} if one can be obtained from the other through the application of a series of the following operations:
\begin{enumerate} 
	\item[(SP):] Symbol permutations: Fix a coordinate position and apply a permutation of the alphabet to all entries in that position.\label{symb}
	\item[(PP):] Positional permutations: Choose two coordinate positions and exchange their entries in every codeword.\label{pos}
\end{enumerate}
If $ C  $ is an $ (n,k,d)_q $ code equivalent to $ C' $, then clearly  $C'$  is also an $ (n,k,d)_q $ code. It follows that any $(n,k)_q$-MDS code is equivalent to an $(n,k)_q$-MDS code which contains the all zero codeword, denoted $ \bar{0} $.

The \textit{Hamming weight}, $ w(c) $ of a codeword, $ c $ is the number of nonzero coordinates it holds (it's distance from  $ \bar{0} $).   For an $(n,k)_q$ code $ C$ we define the \emph{weight spectrum of $ C$} as
$$w( C)=\left\{w(c) \mid c \in  C\setminus \{0\}\right\}.$$

Note that if $ C $ is an $ (n,k)_q $-MDS code containing $ \bar{0} $, then from the MDS property it follows that $ w(C)\subseteq \{n,n-1,\ldots,n-k+1\} $. For a given code $ C $, we denote by $E(w)$  the number of codewords of weight $w$, so that $E(w)=|\{c\in C: wt(c)=w\}|$.   The weight distribution for $(n, k)_q$ MDS codes which contain the zero codeword is completely determined: 
\begin{thm} \cite{395046}
	If $ C $ is an $(n, k)_q$ MDS code which contain the zero codeword, and $ q\ge k $, then:  $E(0)=1$; for $0<w<d$ $E(w)=0$; and for $ w\ge d $
	\begin{equation}
	E(w) = (q-1)\binom{n}{w}\sum_{j=0}^{w-d} (-1)^j\binom{w-1}{j} q^{w-d-j}.\label{weightenumerator}
	\end{equation}
\end{thm}
The weight enumerator can be generalized as follows. Let $T=\{T_1, T_2, \ldots, T_s\}$ be a partition of the set of coordinate positions $\{1, \ldots, n\}$, where $|T_i|=n_i$ for $i=1, \ldots, s$. For each $c\in C$ let $w_i=|supp(c)\cap T_i|$, then $c$ is said to have \emph{$T$-weight profile} 
\[
W_{T}(c)=(w_1, \ldots, w_s).
\] 
The \emph{partition weight enumerator} for a code $C$ with partition $T$ and associated weight profile $W=(w_1, \ldots, w_s)$ is given by 
\[
A^T(W)=\left| \{c\in C\mid W_{T}(c)=(w_1, \ldots, w_s)\}\right|.
\]

The following  is established by El-Khamy and McEliece in  \cite{PWE2005} for linear MDS codes, and in \cite{Alderson2014} for general (not necessarily linear) MDS codes.

\begin{thm}(\cite{PWE2005,Alderson2014})\\ 
	Let $C$ be an $(n, k)_q$ MDS code containing the zero codeword and let $T=\{T_1, \ldots, T_s\}$ be a partition with associated weight profile $W=(w_1, \ldots, w_s)$. Then the partition weight enumerator is given by
	{{\begin{equation}
	A^T(W)=
	(q-1)\binom{n_1}{w_1}\binom{n_2}{w_2}\ldots\binom{n_s}{w_s}\sum_{j=0}^{w-d} (-1)^j\binom{w-1}{j} q^{w-d-j}\label{PWE}
	\end{equation}}}
	where $d\leq w=\sum_{i=1}^s w_i$.
\end{thm}
Both (1) and (2) allow us to compute the weight spectrum of any MDS code, but as observed by Ezermen et. al \cite{Ezerman2011}  it is not immediately obvious  which weights satisfy $ E(w)\ne 0 $. In what follows we shall lay such weight parameters bare.

\section{Some Combinatorial properties of MDS Codes}

Let $C$ be an $(n, k)_q$ MDS code over the alphabet $\A$.   From the MDS property it follows that in any fixed $ k $ coordinate positions, every $ k $-tuple over $ \A $ occurs precisely once as we range over the codewords. In other words, every set of $ k $ coordinates is an \textit{information set}.  It follows that if $ 1\le t\le k $, then for any $ t $  fixed coordinate positions, $ a_1,a_2,\ldots,a_t $ and for any $ \alpha=(\alpha_1,\alpha_2,\ldots,\alpha_t)\in \A^t $  there are precisely $ q^{k-t} $ codewords $ c $ with $ c_{a_i}=\alpha_i $ for all $ i $. Moreover, beyond the $ t $ fixed coordinates, no two such codewords agree in as many as $ k-t $ entries. This property allows us to define $ t $-residual codes:   

\begin{defn}
Let $ C $ be an $ (n,k)_q $-MDS code. Fix $ t \le k $ coordinate positions, $ a_1,a_2,\ldots,a_t $, and   $ \alpha=(\alpha_1,\alpha_2,\ldots,\alpha_t)\in \A^t $. Let $ S $ be the collection of $ q^{k-t} $ codewords $ c $ with $ c_{a_i}=\alpha_i $ for all $ i $.  Deleting the specified $ t $ coordinates from each element of $ S $ results in an $ (n-t,k-t)_q $-MDS code.  Such a code is called a \textit{$ t $-residual code of $ C $}.
\end{defn}  
 
We note that an $ (n,2)_2 $-MDS code $ C $ is easily seen to be combinatorially equivalent to a set of $ n-2 $ mutually orthogonal Latin squares (MOLS), (see e.g. \cite{MR0216968}), or equivalently, a Bruck net of order $ q $ and degree $ n $ \cite{MR0039678}. It follows immediately that $ n\le q+1 $, with equality holding if and only if $ C $ corresponds to an affine plane of order $ q $.      
 
Indeed, MDS codes are rich in combinatorial properties. The following Theorem summarizes some properties we shall find useful in the sequel.

\begin{thm}\label{thm: combinatorial properties}
	Let $C$ be an $(n, k)_q$ MDS code over the alphabet $A$.
	\begin{enumerate}
		\item $ n\ge k $ \label{CP1}
		\item If $ n=k $, then $ C=\A^k $. \label{CP3}
		\item If $ k=1 $ then $ C $ is equivalent to the repetition code. \label{CP4}
		\item If $ C $ contains the zero codeword, then  $ W(C)\subseteq \{n,n-1,\ldots,n-k\} $. \label{CP5}
		\item If $ k>1 $, then $ n\le q+k-1 $. \label{CP6}
		\item If $ q\le k $, then $ n\le k+1 $. \label{CP7}
		\item If $ n<q+k-1 $, and  $\bar{0}\in C $, then  $ n\in W(C) $. \label{CP8}
		\item If $ n=q+k-1 $, $ \bar{0}\in C $, and  $ k,q>2 $, then  $ n\in W(C) $. \label{CP9}
	\end{enumerate}
\end{thm}

\begin{IEEEproof} 
	Parts \ref{CP1}-\ref{CP5}  follow from the definitions and simple counting arguments. \\	Part \ref{CP6} follows by induction on $ k $, with the base case being observed in the previous paragraph. An alternate proof is also given in \cite{MR0216968} (Theorem 4). \\ Part \ref{CP7} follows from a theorem of Bush regarding orthogonal arrays  (\cite{Bush1952}, p.427).  \\ 
Part \ref{CP8}: One could likely perform a careful analysis of (\ref{weightenumerator}) to establish this result, but we prefer a combinatorial argument.   Fix a non-zero $ \alpha \in \A $, and let $ S $ be those codewords with $ \alpha $ in each of the first $ k-1 $ coordinates.  The $ (k-1) $-residual code, $ D $, corresponding to $ S $ is an $ (n-k+1,1)_q $-MDS code, and is therefore equivalent to the repetition code. As such, among all codewords of $ D $ there are $ n-k+1 $ zero entries. Since   $ |D|=q $, and by assumption, $ q>n-k+1 $, the result follows.  \\
Part \ref{CP9}: Again, rather than performing a careful analysis of (\ref{weightenumerator}), we provide a combinatorial argument. By assumption, $ n=q+k-1>k+1 $, so from (\ref{CP7}) we have $ q>k $ and $ n>2k-1=2(k-1)+1 $. According to (\ref{PWE}) there are $ q-1 $ codewords of weight $ d $ that are non-zero in each of the first $ k-1  $ coordinates, let $ c $ be such a codeword. Let $ D $ be the $ (k-1) $-residual code based on the first $ k-1 $ entries of $ c $, where say $ c\mapsto d $. Among all codewords of $ D\setminus \{d\} $ there are $ n-2(k-1)=q-(k-1)<q-1 $ zero entries. So some codeword in $ D $ has full weight. The result follows.

\end{IEEEproof}

\section{Weight Spectra of MDS Codes}

\subsection{The binary MDS codes}

It is well known (see e.g. \cite{MR0465509}) that the only binary linear MDS codes are the trivial ones, namely the $ [n,1]_2 $ (repetition) codes, the  $ [n,n]_2 $ (universe) codes, and the $  [k+1,k]_2 $ single-parity-check codes. We now show that these are in fact all of the binary MDS codes (up to equivalence). Surprisingly, this result does not seem to appear in the literature.

\begin{thm} \label{thm: binary}
If $C$ is an $(n, k)_2$ MDS code, then $ C $ is equivalent to a linear binary MDS code. So $ C $ is equivalent to either the $ [n,1]_2 $ (repetition) codes, the  $ [n,n]_2 $ (universe) code, or the $  [k+1,k]_2 $ single-parity-check code.
 \end{thm}   
\begin{IEEEproof}
 As we are arguing up to equivalence, we may assume with no loss of generality that $ C $ is over the alphabet $ \A=\{0,1\} $, and that $\bar{0}\in C $. \\
 If $ k=1 $, then from Theorem \ref{thm: combinatorial properties} (\ref{CP4}), $ C $ is equivalent to the repetition code. Moreover, since $ \bar{0}\in C $, $ C $ is actually equal to the repetition code $ [n,1]_2 $.\\
  If $ k>1 $, then  from Theorem \ref{thm: combinatorial properties} (\ref{CP1},\ref{CP6}) we have $ k\le n \le k+1 $. We consider these cases separately.\\
 If $ n=k $, then the result follows from Theorem \ref{thm: combinatorial properties} (\ref{CP3}).
For the remaining case, we claim that any $ (k+1,k)_2 $-MDS code  $ C $  is precisely the collection of all words in $ \A^n $ having even weight, so that in fact $ C=[k+1,k]_2 $, the single-parity-check code. Having already  established the case $ k=1 $, we argue inductively. For $ k>1 $ consider the collection $ S $ of codewords with $ 0 $ in the first coordinate. The corresponding residual code  is an $ (n-1,k-1)_2 $-MDS code containing $ \bar{0} $. From the induction hypothesis, all members of $ S $ have even weight. Ranging over all codewords, every binary $ k- $tuple occurs exactly once in the last $ k $ coordinates. All members of $ C\setminus S $ thus have first coordinate $ 1 $, and odd weight in the final $ k $ coordinates. The result follows by induction.     
\end{IEEEproof}

\begin{cor}
	If $ C $ is an $ (n,k)_2 $-MDS code containing $ \bar{0} $ then precisely one of the following hold.
	\begin{enumerate}
		\item $ k=1 $, and $ W(C)=\{n\} $.
		\item $ n=k>1 $, and $ W(C)=\{n,n-1,\ldots,1\} $.
		\item $ n=k+1 $, and $ W(C)=\{1<t\le n \mid t \text{ is even}\}. $
	\end{enumerate}
\end{cor}

\subsection{The 2-Dimensional MDS codes}

The next case we deal with is the 2-dimensional one. As in the linear case, there are two possible weight spectra. 

\begin{thm}\label{thm: 2d} 
	Let $C$ be an $(n, 2)_q$ MDS code containing $ \bar{0} $. 
	\begin{enumerate}
		\item If $ n<q+1 $, then $ W(C)=\{n,n-1\} $.
		\item If $ n=q+1 $, then $ W(C)=\{n-1\} $.		
	\end{enumerate}  
\end{thm}
\begin{IEEEproof}
Let $C$ be an $(n, 2)_q$ MDS code containing $ \bar{0} $. By the MDS property, every nonzero codeword has weight at least $ n-1 $. For each $ i $ with $ 1\le i \le n $, let $ C_i=\{c\in C\mid c\ne \bar{0}, c_i=0\} $.  The $ C_i $'s are mutually disjoint, and $ |C_i|=q-1 $ for each $ i $, so there are precisely $ n(q-1) $ codewords of weight $ n-1 $. The result follows.  	
\end{IEEEproof}

\subsection{The general case for MDS codes}
Having established the base cases, we may now move on to the general setting.  

\begin{thm} \label{thm: general case}
Let  $C$ be an $(n, k)_q$ MDS code containing $ \bar{0} $.   
	\begin{enumerate}
	\item If $ n<q+k-1 $, then $ W(C)=\{n,n-1,\ldots,n-k+1\} $.
	\item If $ n=q+k-1 $, and $ k,q> 2 $ then $ W(C)=\{n, n-1,\ldots,n-k+3,n-k+1\} $.		
\end{enumerate}  
\end{thm}

\begin{IEEEproof}
Let  $C$ be an $(n, k)_q$ MDS code containing $ \bar{0} $. For the first part we argue inductively. The result holds for $ k=1 $ (Theorem \ref{thm: combinatorial properties}, (\ref{CP4})), and for $ k=2 $ (Theorem \ref{thm: 2d}). Fix $ k>2 $, and consider the subset $ S $ of all codewords having first coordinate $ 0 $.  The corresponding residual code $ D $ is an $ (n-1,k-1)_q $-MDS code containing $ \bar{0} $, so $ W(D) (=W(S))=\{n-1,n-2,\ldots,n-k\} $. Theorem \ref{thm: combinatorial properties} (\ref{CP8}) now gives the result.\\
For the second part, let us first consider the case $ k=3 $, so $ C $ is a $ (q+2,3)_q $-MDS code containing $ \bar{0} $. We claim that no codeword has weight $ n-1 $. Indeed, suppose $ c\ne \bar{0}  $ with $ wt(c)<n $. Consider the residual code, $ D $ with respect to a zero coordinate of $ c $, where say $ c\to c' $. From Theorem \ref{thm: 2d} it follows that wt$ (c') ( = wt(c)) = n-2 $. Since $ n\in W(C) $ (Theorem \ref{thm: combinatorial properties} (\ref{CP9})), we have $ W(C)=\{n,n-2\} $. An inductive argument as in the first part now gives the result.          	
\end{IEEEproof}

We would be remiss not to remark that if $ (q+k-1,k)_q $-MDS codes with $k\ge 4 $ do exist, then $ 36 $ divides $ q $  (see \cite{MR704102}). In fact, no such codes are known to exist. 

\begin{cor}
The condition $ (\star) $ holds for an $ (n,k)_q $-MDS code  if $ k=1 $, and holds for $ k>1 $ if and only if $ n<q+k-1 $.
\end{cor}

\subsection{Observations on Distances}

Given an arbitrary $ (n,k)_q $-MDS code, any codeword may be chosen and mapped to $ \bar{0} $ using the operations (SP), or (PP). As such, the weight spectra of non-linear MDS codes translate  to uniform distance properties, or ``distance spectra". To illustrate, Theorems \ref{thm: binary}, \ref{thm: 2d}, and \ref{thm: general case}, and the enumerators (\ref{weightenumerator}), and (\ref{PWE})  give the following.  

\begin{cor}
Let $ C $ be an $ (n,k)_q $-MDS code with $ c\in C $.
\begin{enumerate}
	\item If $ n<q+k-1 $, then there are codewords $ x_0,x_1,\ldots,x_{k-1} \in C$ satisfying $ d(c,x_i)=n-i $ for each $ i $.
	\item If $ n=q+k-1 $, and $ q=2 $ then  there exists a codeword $ x_t \in C$ with $ d(c,x_t)=n-t $ if and only if $ t $ is an even integer with  $  1<t\le n$. 
	\item If $ n=q+k-1 $, and $ q>2 $ then there exists a codeword $ x_t (\ne c) \in C$ with $ d(c,x_t)=n-t $  if and only if  $0\le t\le k-1$, and $ n-t\ne q+1 $.
	
	\item If we denote by $ E(t,c) $ the number of codewords of distance $ t $ from $ c $, then $ E(0,c)=1 $; $ E(t,c)=0 $ for $ 0<t<d $; and for $ t\ge d $ we have
	\[  E(t,c) = (q-1)\binom{n}{t}\sum_{j=0}^{t-d} (-1)^j\binom{t-1}{j} q^{t-d-j}. \] 
	
	\item Let $T=\{T_1, \ldots, T_s\}$ be a partition of the coordinate positions, and let $W=(w_1, \ldots, w_s)$ be an associated weight profile, so $ d\le w=\sum\limits_{i=1}^sw_i\le n $. If $ A^T(W,c) $ denotes the number of codewords that differ from $ c $ in precisely $ w_i $ coordinates among those of $ T_i $ for each $ i $, then
	{{\begin{equation}
			A^T(W,c)=
			(q-1)\binom{n_1}{w_1}\binom{n_2}{w_2}\ldots\binom{n_s}{w_s}\sum_{j=0}^{w-d} (-1)^j\binom{w-1}{j} q^{w-d-j},\label{PWE}
			\end{equation}}}
	where $d\leq w=\sum_{i=1}^s w_i$.
\end{enumerate}  	
	\end{cor}

\section{Conclusion}
In summary we have shown the following:
\begin{itemize}
	\item All binary MDS codes are equivalent to linear.
	\item The condition $ (\star) $ holds for an $ (n,k)_q $-MDS code if $ k=1 $, and for $ k>1 $  if and only if $ n<q+k-1 $.  
	\item An $  (q+k-1,k)_q $-MDS code, $ C $ has the following weight profile 
	\begin{enumerate}
		\item If $ q=2 $, then $ W(C)= W(C)=\{1<t\le n \mid t \text{ is even }\}. $
		\item If $ q>2 $, then $ W(C)=\{n,n-1,\ldots,n-k+1 \}\setminus\{q+1\}. $  
	\end{enumerate} 
\end{itemize}

Moreover, we have the corresponding statements regarding distance spectra.

\section*{Acknowledgment}

The author acknowledges support from the NSERC of Canada through the Discovery Grant Program.

\ifCLASSOPTIONcaptionsoff
  \newpage
\fi



\bibliographystyle{IEEEtran}

\bibliography{IEEEabrv,C://Users/tim/OneDrive/texed/my-texmf/bibtex/bib/main}

%



%
\vspace*{2.5cm}
\end{document}